\documentclass[smallextended]{svjour3}       
\usepackage{graphicx}
\usepackage{amsmath}
\usepackage{amssymb,amsmath,amsfonts}
\usepackage{tabularx}
\usepackage{url}

\usepackage{geometry}
\newcommand{\gph}{\mathrm{gph}\,}
\newcommand{\Limsup}{\mathop{{\rm Lim}\,{\rm sup}}}
\newcommand{\dom}{\mbox{\rm dom}\,}
\newcommand{\argmin}{\mathop{\rm arg\,min}\limits}
\newcommand{\inn}{{\rm int\,}}
\newcommand{\boxatend}{\hfill$\square$}
\begin{document}

\title{On computation of optimal strategies in oligopolistic markets respecting the
cost of change\thanks{The research of J. V. Outrata was supported by the Czech Science Foundation (GA ČR), through the grant 17-08182S and by the Australian Research Council, project DP160100854F. The research of J. Valdman  was supported by the Czech Science Foundation (GA ČR), through the grant 17-04301S.
}
}


\author{Ji\v{r}\'{i} V. Outrata         \and
        Jan Valdman 
}


\institute{Ji\v{r}\'{i} V. Outrata \at
              Department of Decision-Making Theory, Czech Academy of Sciences, Institute of Information Theory and Automation, Prague, Czech Republic  \& Federation University of Australia, School of Science, Information Technology and Engineering, Ballarat, Australia \\
              \email{outrata@utia.cas.cz}           
           \and
           Jan Valdman \at
             Department of Decision-Making Theory, Czech Academy of Sciences, Institute of Information Theory and Automation, Prague 
\&
Institute of Mathematics, Faculty of Science, University of South Bohemia, \v{C}esk\'e Bud\v{e}jovice, Czech Republic \\
\email{jan.valdman@utia.cas.cz}
}

\date{Received: date / Accepted: date}

\maketitle

\begin{abstract}
The paper deals with a class of parametrized equilibrium problems, where the objectives of the players do possess nonsmooth terms. The respective Nash equilibria can be characterized via a parameter-dependent variational inequality of the second kind, whose Lipschitzian stability is thoroughly investigated. This theory is then applied to evolution of a oligopolistic market in which the firms adopt their production strategies to changing input costs, while each change of the production is associated with some "costs of change". We examine both the Cournot-Nash equilibria as well as the two-level case, when one firm decides to take over the role of the Leader (Stackelberg equilibrium). The impact of costs of change is illustrated by academic examples.
\end{abstract}

\section{Introduction}
Consider an oligopolistic market, where the data of the production cost functions and constraints
of each producer/firm are available to all his rivals. In such a case each producer can compute his
optimal non-cooperative Cournot-Nash strategy by solving the corresponding variational inequality,
see \cite{MSS} and \cite{OKZ}. It may happen, however, that in the course of time some external
parameters change, e.g., the prices of the inputs or the parameters of the inverse demand function
describing the behavior of the customers. In such a case, the strategies should be adjusted but, as
thoroughly analyzed in \cite{Fl}, each change is generally  associated with some expenses, called
{\em costs of change.} Thus, given a certain {\em initial strategy profile} (productions of all
firms), we face then a different equilibrium model, in which the costs of change enter the
objectives of some (or all) producers. Since these costs are typically nonsmooth, the
respective variational inequality, describing the new non-cooperative equilibrium, becomes
substantially more complicated, both from the theoretical as well as from the numerical point of
view. One can imagine that such updates of strategies are performed repetitively. This leads to a
discrete-time evolution process, where the firms respond to changing conditions by repetitive
solution of the mentioned rather complicated variational inequality (with updated data). As
discussed in \cite[Chapter 12]{OKZ}, it may also happen that one of the producers, having an advantage over the others, takes over  the role of a Leader and switches to the {\em Stackelberg
strategy}, whereas the remaining firms continue to play non-cooperatively with each other on the lower level as
Followers. In this case, our discrete-time evolution process, considered from the point of view of
the Leader, amounts to repetitive solution of {\em mathematical programs with equilibrium
constraints} (MPECs) with a nonsmooth objective and the above mentioned variational inequality
among the constraints.
Further, it is interesting to note that the above described model has, in case of  positively homogeneous costs of change, a similar structure as 
some infinite-dimensional variational systems used in continuum mechanics to model a class of {\em rate-independent processes} cf., e.g., \cite{MR} or \cite{FKV}.

The plan of the paper is as follows. In the preliminary Section 2 we collect the necessary background from variational analysis. Section 3 consists of two parts. In the first one we introduce a general parameter-dependent non-cooperative equilibrium problem which is later used for modeling of the considered oligopolistic market. By employing standard arguments, existence of the respective solutions (equilibria) is shown. In the second part we then consider a parameter-dependent variational system which encompasses the mentioned equilibrium problem and is amenable to advanced tools of variational analysis. In this way one obtains useful stability results concerning the respective {\em solution map}. They are used in the sequel but they are also important for its own sake.

Thereafter, in Section 4, this equilibrium problem is specialized to a
form, corresponding to the oligopolistic market model from \cite{MSS}.
In this case, the solution map is even single-valued and locally Lipschitzian. In Section 5 we then consider a modification of the 5-firm example from \cite{MSS} with the aim to illustrate the role of costs of change and to describe a possible numerical approach to the computation of the respective equilibria. Whereas Section 5.1 deals with noncooperative Cournot-Nash equilibrium, Section 5.2 concerns
 the situation, when one of the firms prefers to apply the Stackelberg strategy.
 In both cases our main numerical tool is the so-called Gauss-Seidel method from \cite{Ka} tailored to the considered type of problems.\\

The following notation is employed. Given a vector $x \in \mathbb{R}^{n}, [x]$ denotes the linear subspace generated by $x$ and $[x]^{\perp}$stands for its orthogonal complement. For a multifunction $F: \mathbb{R}^{n} \rightrightarrows \mathbb{R}^{m}, \gph F$ signifies the graph of $F$, $\delta_{A}$ is the indicatory function of a set $A$ and $\overline{\mathbb{R}} = R \cup \{+ \infty \}$ is the extended real line. $\mathbb{B}$ stands for the unit ball and, for a cone $K$, $K^{\circ}$  denotes its (negative) polar.

 \section{Background from variational analysis}
Throughout the whole paper, we will make an extensive use of the following basic notions of modern variational analysis.
\begin{definition}\label{Def.1}
Let $A$ be a closed set in $\mathbb{R}^{n}$ and $\bar{x}\in A$. Then
\[
T_{A}(\bar{x}):=\Limsup\limits_{t \searrow 0}    \frac{A-\bar{x}}{t}
\]
is the {\em tangent (contingent, Bouligand)} cone to $A$ at $\bar{x}$,
$$\widehat{N}_{A}(\bar{x}):= (T_{A}(\bar{x}))^{\circ}
$$
is {\em regular (Fr\'{e}chet) normal cone} to $A$ at $\bar{x}$, and
$$
N_{A}(\bar{x}):= \Limsup\limits_{A \atop x\rightarrow \bar{x}} \widehat{N}_{A}(x) =
\{x^{*}|\exists \stackrel{A}{x^{*}_{i} \rightarrow x^{*}}, x^{*}_{i}\in \widehat{N}_{A}(x_{i}) \mbox{ such that } x^{*}_{i} \rightarrow x^{*}\}
$$
 is the {\em limiting (Mordukhovich) normal cone} to $A$ at $\bar{x}$.
\end{definition}
In this definition "$\Limsup$" stands for the  Painlev\'{e}-Kuratowski {\em outer set limit}.
If $A$ is convex, then $\widehat{N}_{A}(\bar{x})=N_{A}(\bar{x})$ amounts to the  classical normal cone in the sense of convex analysis and we  write $N_{A}(\bar{x})$.

In the sequel, we will also employ the co-called critical cone. In the setting of Definition \ref{Def.1} with an given normal $d^{*}\in \widehat{N}_{A}(\bar{x})$, the cone
\[
\mathcal{K}_{A}(\bar{x},d^{*}):= T_{A}(\bar{x})\cap [d^{*}]^{\perp}
\]
is called the {\em critical cone} to $A$ at $\bar{x}$ with respect to $d^{*}$.

The above listed cones enable us to describe the local behavior of set-valued maps via various generalized derivatives. Consider a closed-graph multifunction $F$ and the point $(\bar{x},\bar{y})\in \gph F$.
\begin{definition}\label{Def.2}
 \begin{enumerate}
 \item [(i)] The multifunction $DF (\bar{x},\bar{y}):\mathbb{R}^{n}   \rightrightarrows \mathbb{R}^{m}$, defined by
     \[
     DF(\bar{x},\bar{y})(u):=\{v \in \mathbb{R}^{m} | (u,v)\in T_{\gph F}(\bar{x},\bar{y})\}, d \in \mathbb{R}^{n},
     \]
     is called the {\em graphical derivative} of $F$ at $(\bar{x},\bar{y})$;
 \item [(ii)] The multifunction $D^{*}F(\bar{x},\bar{y}): \mathbb{R}^{m}   \rightrightarrows \mathbb{R}^{n}$, defined by
 \[
      D^{*}F(\bar{x},\bar{y})(v^{*}):=\{u^{*} \in \mathbb{R}^{n} | (u^{*},-v^{*})\in N_{\gph F}(\bar{x},\bar{y})\}, v^{*} \in \mathbb{R}^{m},
     \]
     is called the {\em limiting (Mordukhovich) coderivative} of $F$ at $(\bar{x},\bar{y})$.

\end{enumerate}
\end{definition}

Next we turn our attention to a proper convex, lower-semicontinuous  (lsc) function $q:\mathbb{R}^{n}\rightarrow\overline{\mathbb{R}}$. Given an $\bar{x} \in \dom q$, by $\partial q(\bar{x})$ we  denote the classical Moreau-Rockafellar {\em subdifferential} of $q$ at $\bar{x}$. In this case, for the {\em subderivative} function $dq(\bar{x}):\mathbb{R}^{n}\rightarrow\overline{\mathbb{R}}$ (\cite[Definition 8.1]{RW}) it holds that
\[
dq(\bar{x})(w)=  q^{\prime}(\bar{x};w) :=  \lim\limits_{\tau\searrow 0}\frac{q(\bar{x}+\tau w)-q(\bar{x})}{\tau}
 \mbox{ for all }~~w \in \mathbb{R}^{n}.
\]
  In Section 3 we will employ also second-order subdifferentials and second-order subderivatives of $q$.
\begin{definition}\label{Def.3}
Let $\bar{v}\in \partial q(\bar{x})$. The multifunction $\partial^{2}q(\bar{x},\bar{v}):\mathbb{R}^{n}  \rightrightarrows\mathbb{R}^{n}$ defined by
\[
\partial^{2}q(\bar{x},\bar{v})(v^{*}):= D^{*}\partial q(\bar{x},\bar{v})(v^{*}), ~~ v^{*}\in \mathbb{R}^{n},
\]
 is called the {\em second-order subdifferential} of $q$ at $(\bar{x},\bar{v})$.
\end{definition}

If $q$ is separable, i.e., $q(x)=\sum\limits^{n}_{i=1}q_{i}(x_{i})$ with some proper convex, lsc functions $q_{i} : \mathbb{R}\rightarrow\overline{\mathbb{R}} $, then
\[
 \partial^{2} q(\bar{x},\bar{v})(v^{*}) =
 \left [ \begin{array}{c}
 \partial^{2} q_{1}(\bar{x}_{1},\bar{v}_{1})(v^{*}_{1})\\
 \vdots \\
  \partial^{2} q_{n}(\bar{x}_{n},\bar{v}_{n})(v^{*}_{n})
 \end{array}
\right ],
\]
where $\bar{v}_{i}, v^{*}_{i}$ are the {\em i}th components of the vectors, $\bar{v},v^{*}$, respectively. \\

Concerning second-order subderivatives (\cite[Definition 13.3]{RW}), we confine ourselves to the case when $q$  is, in addition, {\em piecewise linear-quadratic}. This means that $\dom q$ can be represented as the union of finitely many polyhedral sets, relative to each of which $q(x)$ is given in the form $\frac{1}{2}\langle x,Ax \rangle + \langle a,x \rangle  + \alpha $ for some scalar $\alpha \in \mathbb{R}$, vector $a \in \mathbb{R}^{n}$ and a symmetric $[n \times n]$ matrix $A$, cf. \cite[Definition 10.20]{RW}.

In this particular case it has been proved in   \cite[Proposition 13.9]{RW} that, with $\bar{v}\in \partial q(\bar{x})$ and $w \in \mathbb{R}^{n}$ the second-order subderivative $d^{2}q(\bar{x} | \bar{v})$ is proper convex and piecewise linear quadratic and
\begin{equation}\label{eq-102}
d^{2}q(\bar{x} | \bar{v})(w)= q^{\prime \prime} (\bar{x};w)+ \delta_{K(\bar{x},\bar{v})}(w),
\end{equation}
where
\[
q^{\prime \prime} (\bar{x};w) := \lim\limits_{\tau\searrow 0}
\frac{q(\bar{x}+\tau w)-q(\bar{x}) - \tau q^{\prime}(\bar{x};w)}
{\frac{1}{2}\tau^{2}}
\]
is the {\em one-sided second directional derivative} of $q$ at $\bar{x}$ in direction $w$
and $K(\bar{x},\bar{v}):=\{w | q^{\prime}(\bar{x};w)= \langle \bar{v},w\rangle \}$. For a general theory of second-order subderivatives (without our restrictive requirements) the interested reader is referred to \cite[Chapter  13 B]{RW}.

We conclude now this section with the definitions of two important Lipschitzian stability notions for multifunctions which will be extensively employed in the sequel.
\begin{definition}\label{Def.4}
Consider a multifunction $S:\mathbb{R}^{m} \rightrightarrows \mathbb{R}^{n}$ and a point $(\bar{u},\bar{v})\in \gph S$.
\begin{enumerate}
 \item [(i)]
 $S$ is said to have the {\em Aubin property} around $(\bar{u},\bar{v})$, provided there are neighborhoods $\mathcal{U}$ of $\bar{u}$, $\mathcal{V}$ of $\bar{v}$ along with a constant $\eta \geq 0$ such that
 \[
 S(u_{1}) \cap \mathcal{V} \subset S(u_{2})+\eta \| u_{1}-u_{2}\| \mathbb{B} \mbox{ for all } u_{1},u_{2}\in \mathcal{U}.
 \]
 \item [(ii)]
 We say that  $S$ has a
 {\em single-valued and Lipschitzian localization} around $(\bar{u},\bar{v})$, provided there are neighborhoods $\mathcal{U}$ of $\bar{u}$, $\mathcal{V}$ of $\bar{v}$ and a Lipschitzian mapping $s:\mathcal{U} \rightarrow \mathbb{R}^{n}$ such that $s(\bar{u})=\bar{v}$ and
 \[
 S(u)\cap \mathcal{V} = \{s (u)\} \mbox{ for all } u \in \mathcal{U}.
 \]
\end{enumerate}
\end{definition}

\noindent Further important stability notions can be found, e.g., in \cite{DR}.

\section{General equilibrium model: Existence and stability}
Consider a non-cooperative game of $l$ players, each of which solves the optimization problem 
\begin{equation}\label{eq-1}
\begin{array}{ll}
\mbox{ minimize } & f_{i}(p,x_{i},x_{-i}) + q_{i}(x_{i})\\
\mbox{ subject to } & \\
& x_{i} \in A_{i},
\end{array}
\end{equation}
$i = 1,2,\ldots, l$. In (\ref{eq-1}), $x_{i}\in \mathbb{R}^{n}$ is the {\em strategy} of the
$i$th player,
$$x_{-i} := (x_{1}, \ldots, x_{i-1}, x_{i+1}, \ldots, x_{l})
\newline  \in (\mathbb{R}^{n})^{l-1}$$ is the
{\em strategy profile} of the remaining players and $p \in \mathbb{R}^{m}$  is a parameter, common
for all players. Further, the functions
$$f_{i}:\mathbb{R}^{m}\times (\mathbb{R}^{n})^{l}\rightarrow
\mathbb{R} \quad
\mbox{ and }
\quad q_{i}: \mathbb{R}^{n}\rightarrow \mathbb{R}, \quad i = 1,2,\ldots, l,$$  are continuously
differentiable and convex continuous, respectively, and the sets of {\em admissible strategies}
$A_{i}, i=1,2,\ldots, l$, are closed and convex. The objective in (\ref{eq-1}) is thus the sum of a
smooth function depending on the whole strategy profile $x:=(x_{1},x_{2},\ldots, x_{l})$ and a
convex (not necessarily smooth) function depending only on $x_{i}$. Let us recall that, given a
parameter vector $\bar{p}$, the strategy profile $\bar{x}=(\bar{x}_{1}, \bar{x}_{2}, \ldots,
\bar{x}_{l})$ is a corresponding {\em Nash equilibrium} provided 
\[
\bar{x}_{i} \in \argmin_{x_{i}\in A_{i}}
\big[ f_{i}(\bar{p},x_{i},\bar{x}_{-i})+q_{i}(x_{i}) \big] \mbox{ for all } i.
\]
Denote by $S:\mathbb{R}^{m}\rightrightarrows (\mathbb{R}^{n})^{l}$ the solution mapping which assigns
each $p$ the corresponding  (possibly empty) set of Nash equilibria.
The famous Nash Theorem \cite[Theorem 12.2]{Au} yields the next statement.
 \begin{theorem}\label{Thm.1}
 Given $\bar{p} \in \mathbb{R}^{n}$, assume that
 \begin{enumerate}
 \item [(A1)]
 for all admissible values of $x_{-i}$
 functions $f_{i}(\bar{p}, \cdot, x_{-i}), i=1,2,\ldots, l$, are convex,   and
 \item [(A2)]
 sets $A_{i}, i=1,2,\ldots, l$, are bounded.
\end{enumerate}
Then  $S(\bar{p})\neq \emptyset$.

 \end{theorem}
Suppose from now on that (A1) holds true for all $p$ from an open set $\mathcal{B}\subset \mathbb{R}^{m}$.
Then one has that $\mathcal{B}\subset \dom S$ and for $p \in \mathcal{B}$ 
\begin{equation}\label{eq-2}
S(p)=\{x | \, 0 \in F(p,x)+Q(x)\},
\end{equation}
 where
\[
\begin{split}
& F(p,x)=  \left[ \begin{array}{c}
F_{1}(p,x) \\
\vdots \\
F_{l}(p,x)
\end{array}\right ]
\mbox{ with } F_{i}(p,x)=\nabla_{x_{i}}f_{i}(p,x_{i},x_{-i}), \, i=1,2,\ldots, l, \,\,\mbox{ and
}\\[2ex]
& Q(x)=  \partial \tilde{q}(x) \mbox{ with } \tilde{q}(x)=\sum\limits^{l}_{i=1} \tilde{q}_{i}(x_{i}) \mbox{ and } \tilde{q}_{i}(x_{i})= q_{i}(x_{i})+\delta_{A}(x_{i}), \,  i = 1,2,\ldots, l.
\end{split}
\]
This follows immediately from the fact that under the posed assumptions the solution set of
(\ref{eq-1}) is characterized by the first-order condition 
\[
0 \in \nabla_{x_{i}}f_{i}(p, x_{i},x_{-i}) + \partial \tilde{q}_{i}(x_{i}), ~ i = 1,2,\ldots, l.
\]
\if{
 Next we state two results concerning the local stability and sensitivity  of $S$ around a given
reference point. They are not indispensable for the numerical technique developed in Section 4, but
the first one will be helpful in the ImP approach of Section 5.

Consider a parameter $\bar{p}\in \mathcal{B}$ and an $\bar{x} \in S(\bar{p})$. Further, put
$\bar{v}:= -F(\bar{p},\bar{x})$ and $\bar{v}_{i}= -F_{i}(\bar{p},\bar{x}), ~ i=1,2,\ldots, l$.

\begin{theorem}\label{Thm.1}
Assume that the {\em adjoint} GE 
\begin{equation}\label{eq-3}
0 \in \sum\limits^{l}_{i=1}\nabla_{x}F_{i}(\bar{p},\bar{x})^{T}u_{i}+
\left[ \begin{array}{c}
\partial^{2}\tilde{q}_{1}(\bar{x}_{1},\bar{v}_{1})(u_{1}) \\
\vdots\\
 \partial^{2}\tilde{q}_{l}(\bar{x}_{l},\bar{v}_{l})(u_{l})
\end{array}\right]
\end{equation}
in variable $u=(u_{1},u_{2}, \ldots, u_{l})\in (\mathbb{R}^{n})^{l}$ has only the trivial
solution $u=0$. Then $S$ has the Aubin property around $(\bar{p},\bar{x})$.
\end{theorem}

\proof
The statement follows immediately from \cite[Section 4.4]{M1} taking into account that, thanks to
the separability of $Q$, 
\[
D^{*}Q(\bar{x},\bar{v})(u)=
\left[ \begin{array}{c}
D^{*}\partial\tilde{q}_{1}(\bar{x}_{1},\bar{v}_{1})(u_{1}) \\
\vdots\\
D^{*} \partial\tilde{q}_{l}(\bar{x}_{l},\bar{v}_{l})(u_{l})
\end{array}\right] =
\left[ \begin{array}{c}
\partial^{2}\tilde{q}_{1}(\bar{x}_{1},\bar{v}_{1})(u_{1}) \\
\vdots\\
 \partial^{2}\tilde{q}_{l}(\bar{x}_{l},\bar{v}_{l})(u_{l})
\end{array}\right].
\]
\endproof

In applications one has thus to compute the second-order subdifferentials of functions
$\tilde{q}_{i}$. If $n=1$, this can be done on the basis of the results from \cite{MorOut01}.
Denote by $\Sigma$ the mapping which assigns a pair $(\xi,p)\in (\mathbb{R}^{n})^{l}\times
\mathbb{R}^{m}$ the set of solutions to the GE 
\[
\xi \in F(p,x)+Q(x).
\]
Then the condition of Theorem \ref{Thm.1} is necessary and sufficient for the Aubin property of
$\Sigma$ around $(0,\bar{p},\bar{x})\in \gph \Sigma$.

 \begin{theorem}\label{Thm.2}
Let the assumption of Theorem \ref{Thm.1}  be fulfilled. Further assume that the functions $q_{i},
i=1,2,\ldots, l$, are piecewise linear quadratic, cf. \cite[Definition 10.20]{RW}. Then for all
directions $h \in \mathbb{R}^{m}$ one has 
\begin{equation}\label{eq-4}
DS(\bar{p},\bar{x})(h)=\left\{k ~|~ \nabla_{p}F(\bar{p},\bar{x})h + \nabla_{x}F(\bar{p},\bar{x})k
+
 \left[ \begin{array}{c}
\partial \varphi_{1}(k_{1}) \\
\vdots \\
\partial \varphi_{l}(k_{l})
\end{array}\right ]
\right\},
\end{equation}

where $\varphi_{i}(k)= \frac{1}{2}d^{2}\tilde{q}_{i} (\bar{x}_{i}, \bar{v}_{i})(k), ~
i=1,2,\ldots, l$.
\end{theorem}

\proof
In terms of GE (\ref{eq-2}) the condition imposed in Theorem \ref{Thm.1} implies the fulfillment
of the standard qualification condition (cp. \cite[Theorem 6.14]{RW})
\begin{equation}\label{eq-5}
\left. \begin{split}
& 0= \nabla_{p}F(\bar{p},\bar{x})^{T}u\\ & 0 \in \nabla_{x}F(\bar{p},\bar{x})^{T}u +
D^{*}Q(\bar{x},\bar{v})(u)
\end{split}\right\}\Rightarrow u =0.
\end{equation}
Clearly,
\begin{equation}\label{eq-201}
\gph S = \left\{(p,x) \left|
\left[ \begin{split}
& x \\ - & F(p,x)
\end{split}\right ] \in \gph Q  \right. \right\}
\end{equation}
and implication (\ref{eq-5}) ensures that the constraint system on the right-hand side of
(\ref{eq-201}) fulfills the  Abadie constraint qualification at $(\bar{p},\bar{x})$, see
\cite[Proposition 1]{HO}. It follows that 
\[
T_{\gph S}(\bar{p},\bar{x})= \left\{(h,k) \left|
\left[ \begin{split}
& k \\ - & \nabla F(\bar{p},\bar{x})(h,k)
\end{split}\right ] \in T_{\gph Q}(\bar{x}, -F(\bar{p},\bar{x}))
 \right.\right\}
\]
 or, equivalently, 
\[
DS(\bar{p},\bar{x})(h)= \{k | 0 \in \nabla f(\bar{p},\bar{x})(h,k)+ DQ(\bar{x},\bar{u})(k)\}
\]
for all $h \in \mathbb{R}^{s}$. So it remains just to compute $DQ(\bar{x},\bar{u})(k)$. To this
aim observe first that, thanks to the assumptions imposed on $q_{i}$ and $A_{i}$, functions
$\tilde{q}_{i}$ are {\em fully amenable} \cite[Example 10.24]{RW}. This implies that they are {\em
twice epi-differentiable, prox-regular} and {\em subdifferentially continuous} \cite[Theorem 13.14
and Proposition 13.32]{RW} and so we may invoke \cite[Theorem 13.40]{RW}, according to which the
sets $\gph \partial \tilde{q}_{i}, i = 1,2,\ldots, l$, are {\em geometrically derivable}. Since
$\gph Q = \XX\limits^{l}_{i=1} \gph \partial \tilde{q}_{i}$, the inclusion in \cite[Proposition
6.41]{RW} becomes equality 
\[
T_{\gph Q}(\bar{x},\bar{v})= \XX\limits^{l}_{i=1}T_{\gph}\partial \tilde{q}_{i}
(\bar{x}_{i},\bar{v}_{i}),
\]
which implies in turn that 
\begin{equation}\label{eq-14}
DQ(\bar{x},\bar{u})(k)=
\left[ \begin{array}{c}
D\partial\tilde{q}_{1}(\bar{x}_{1},\bar{v}_{1})(k_{1}) \\
\vdots\\
D \partial\tilde{q}_{l}(\bar{x}_{l},\bar{v}_{l})(k_{l})
\end{array}\right], ~ k \in \mathbb{R}^{s}.
\end{equation}
Next we recall \cite[Theorem 13.40]{RW} once more and conclude that the graphical derivatives on
the right-hand side of (\ref{eq-14}) can be expressed in terms of second subderivatives via  the
relations 
\begin{equation}\label{eq-15}
D \partial\tilde{q}_{i}(\bar{x}_{i},\bar{v}_{i}) = \partial \varphi_{i} \mbox{ for } \varphi_{i} =
\frac{1}{2}d^{2}\tilde{q}_{i} (\bar{x}_{i}, \bar{v}_{i}), ~ i=1,2,\ldots, l.
\end{equation}
The proof is complete.
\endproof

The satisfaction of the assumption of Theorem \ref{Thm.1} can be sometimes ensured in a very simple
way. 
\begin{corollary}\label{Cor.1}
Assume in the setting of Theorem \ref{Thm.1} that the $[nl \times nl]$ matrix
$\nabla_{x}F(\bar{p},\bar{x})$ is positive definite. Then $S$ has the Aubin property around
$(\bar{p},\bar{x})$. Moreover, if $q_{i}$   are piecewise linear quadratic, $i-1,2,\ldots, l$, then
$\bar{x}$ is locally unique in $S(\bar{p})$.
\end{corollary}

\proof
Let us premultiply GE (\ref{eq-3}) from the left-hand side by $u=(u_{1}, u_{2}, \ldots, u_{l})$
which yields the relation
\begin{equation}\label{eq-16}
0 \in \langle \nabla_{x}F(\bar{p},\bar{x})u,u \rangle + \sum\limits^{l}_{i=1}
\langle u_{i}, \partial^{2}\tilde{q}_{i}(\bar{x}_{i},\bar{v}_{i}) (u_{i})\rangle .
\end{equation}
 The second term on the right-hand side of (\ref{eq-16}) amounts to the product $\langle u,D^{*}
\partial\tilde{q}(\bar{x},\bar{v}) u\rangle $, where $\tilde{q}(x)= \sum\limits^{l}_{i=1}
\tilde{q}_{i}
(\bar{x}_{i})$. Since $\tilde{q}$ is convex by the assumptions, the mapping $\partial\tilde{q}:
\mathbb{R}^{nl} \rightrightarrows \mathbb{R}^{nl}$ is maximal monotone \cite[Theorem 12.17]{RW}. We
can thus invoke \cite[Theorem 2.1]{PR} according to which this term is nonnegative. It follows that
(\ref{eq-16}) possesses only the trivial solution $u=0$, which proves the first assertion. To prove
the second one, observe that, thanks to \cite[Theorem 13.57]{RW}, one has 
\[
DQ(\bar{x},\bar{v})(u)\subset D^{*}Q(\bar{x},\bar{v})(u)
\]
 for all $u \in (\mathbb{R}^{n})^{l}$. This implies by the same argument as above that the
graphical derivative 
\[
DS(\bar{p},\bar{x})(0)= \{k | 0 \in \nabla_{x}F(\bar{p},\bar{x})k+DQ(\bar{x},\bar{v})(k)\}
\]
 has only the trivial solution $k = 0$. From the Levy-Rockafellar criterion \cite[Theorem
4C.1]{DR} it follows now that $S$ possesses also the isolated calmness property at
$(\bar{p},\bar{x})$ which entails the local uniqueness of $\bar{x}$ in $S(\bar{p})$. The statement
is proved.
\endproof

Let us illustrate the results of Theorem \ref{Thm.1} and Corollary \ref{Cor.1} via a simple
academic example.

\begin{example}
Put $s=2, n=l=1$ and consider the GE 
\[
0 \in p_{1} + p_{2}x + \partial \tilde{q}(x),
\]
where $\tilde{q}(x)= |x| + \delta_{A}(x)$ with $A = [0,1]$. Consider the reference pair
$((\bar{p}_{1}, \bar{p}_{2}),\bar{x})= ((-1,1),0)$ so that $\bar{v}=1$. Since
$\nabla_{x}F(\bar{p},\bar{x})=\bar{p}_{2}=1$, Corollary \ref{Cor.1} can be applied and we conclude
that $S$ has the Aubin property around $(\bar{p},\bar{x})$. Moreover, the  absolute value function
is piecewise linear quadratic,  and so $\bar{x}$ is locally unique in $\bar{p}$. To compute
$DS(\bar{p},\bar{x})$ we recall \cite[Definition 13.3]{RW} according to which
\[
\varphi(k):=d^{2}q(\bar{x},\bar{v})(k)= \delta_{\mathbb{R}_{+}}(k),
\]
 and thus
\begin{equation}\label{eq-17}
DS(\bar{p},\bar{x})(h_{1},h_{2})=\{k | 0 \in h_{1}+k+ N_{\mathbb{R}_{+}}(k)\}.
\end{equation}
The GE in (\ref{eq-17}) amounts to the complementarity problem 
\[
k \geq 0, ~ h_{1}+k \geq 0, ~ \langle h_{1}+k,k\rangle =0
\]
 and we may conclude that 
\[
DS(\bar{p},\bar{x})(h_{1},h_{2})=
\left \{ \begin{array}{cl}
-  h_{1} & \mbox{ if } h_{1}< 0\\ 0 & \mbox{ otherwise },
\end{array}\right.
\]
 see\\
  Let us illustrate the results of Theorem \ref{Thm.1} and Corollary \ref{Cor.1} via a simple academic example.
\begin{figure}[h]
\centering
\includegraphics[width=0.48\textwidth]{figure_S}
\includegraphics[width=0.48\textwidth]{figure_DS}
\caption{The set of equilibria $S(p_1, p_2)$ (left) and the graphical derivative $DS(h_1,h_2)$ evaluated at the point $p_1=-1, p_2=1$ of  Example 1.}
\label{mesh}
\end{figure}

\end{example}
On the basis of  Corollary \ref{Cor.1} we can now prove the following non-local result. 
 \begin{theorem}\label{Thm.3}
Assume that (A1) is fulfilled for all $p$  from some open set $\mathcal{B}$, (A2) holds true and
$\nabla_{x}F(p,x)$ is positive definite over $\mathcal{B} \times \XX\limits^{l}_{i=1} A_{i}$. Then
$S$ is single-valued and locally Lipschitz on $\mathcal{B}$.
  \end{theorem}
  \proof
 As already  mentioned, under the posed assumptions, $S$  is nonempty-valued on $\mathcal{B}$. Next
 we observe that, by virtue of \cite[Theorem 3.43]{OR}, for  each $p \in \mathcal{B}$ the operator
 $F(p,\cdot)$ is strictly monotone over $\XX\limits^{l}_{i=1} A_{i}$. Since the multifunction $Q$
 amounts to the subdifferential of the proper, convex function $\tilde{q}_{1}(x_{1})+ \ldots +
 \tilde{q}_{l}(x_{l}), Q$ is monotone (\cite[Theorem 12.17]{RW}). It follows by \cite[Exercise
 12.4c]{RW} that for each $p \in \mathcal{B}$ the multifunction
 \[
 G(p, \cdot):=F(p, \cdot) + Q(\cdot)
 \]
 is strictly  monotone over   $\XX\limits^{l}_{i=1} A_{i}$. This enables us to invoke a complement
 to
 \cite[Theorem 12.51]{RW} (mentioned at the end of Section 12G) and state that $S$ is in fact
 single-valued on $\mathcal{B}$.

 From Corollary \ref{Cor.1} we know that $S$ has the Aubin property around all points $(p,x)\in
 \gph S \cap (\mathcal{B} \times (\mathbb{R}^{n})^{l})$. Since this property amounts to the local
 Lipschitz continuity for single-valued maps, the statement has been established.
  \endproof
\fi{

  \begin{remark}\label{Rem.1}
  The GE in (\ref{eq-2}) can be equivalently written down in the form:

  For a given $\bar{p}$ find $\bar{x}$ such that
  \[
  \langle F(\bar{p},\bar{x}), x - \bar{x}\rangle + \tilde{q}(x)-\tilde{q}(\bar{x})\geq 0 ~ {\rm  for ~all} ~ x.
  \]

  Our equilibrium is thus governed by a variational inequality (VI) of the second kind, cf. \cite{KS}, \cite[page 96]{FP}.
  \end{remark}

  Next we will concentrate on the (local) analysis of $S$ (given by (\ref{eq-2})) under the less restrictive assumptions that, with $s:= ln$,
  \begin{description}
  \item [(i)]
  $F:\mathbb{R}^{m}\times\mathbb{R}^{s}\rightarrow\mathbb{R}^{s}$ is continuously differentiable, and
  \item [(ii)]
  $Q(\cdot)= \partial \tilde{q} (\cdot)$ for a proper convex, lsc function $\tilde{q}:\mathbb{R}^{s} \rightarrow \overline{\mathbb{R}}$.
  \end{description}
  In this way, the obtained results will be applicable not only to the equilibrium problem stated above, but to a broader class of parametrized VIs of the second kind. Note that stability of the {\em generalized equation} (GE)
  \begin{equation}\label{eq-50}
0 \in F(p,x)+\partial \tilde{q}(x)
  \end{equation}
  has been investigated, among other works, in \cite[Chapter 13]{RW} even without any convexity assumptions imposed on  $\tilde{q}$. As proved in \cite[Theorem 13.48]{RW}, $S$ has the Aubin property around $(\bar{p},\bar{x}) \in \gph S$ provided the {\em adjoint} GE
  \begin{equation}\label{eq-51}
 0 \in \nabla_{x}F(\bar{p},\bar{x})^{T}u + \partial^{2}\tilde{q}(\bar{x},-F(\bar{p},\bar{x}))(u)
  \end{equation}
in variable $u \in \mathbb{R}^{s}$ has only the trivial solution $u = 0$.

This condition is automatically fulfilled provided $\nabla_{x}F(\bar{p},\bar{x})$ is positive definite. Indeed, due to the assumptions imposed on $\tilde{q}$, the mapping $\partial \tilde{q}$ is maximal monotone \cite[Theorem 12.17]{RW}. We can thus invoke \cite[Theorem 2.1]{PR}, according to which
\[
\langle u,v \rangle \geq 0~ \mbox{ for all } ~v \in \partial^{2}\tilde{q} (\bar{x},-F(\bar{p},\bar{x}))(u).
\]
The result thus follows from the inequality
\[
0 \leq \langle u,\nabla_{x}F(\bar{p},\bar{x})^{T}u\rangle.
\]
Let us now derive conditions ensuring the existence of a single-valued and Lipschitzian localization of $S$ around $(\bar{p},\bar{x})$. To this purpose we employ \cite[Theorem 3G.4]{DR}, according to which this property of $S$ is implied by the existence of a single-valued and Lipschitzian localization of the associated partially linearized mapping $\Sigma : \mathbb{R}^{s}\rightrightarrows  \mathbb{R}^{s}$ defined by
\begin{equation}\label{eq-52}
\Sigma(w):=\left \{ x | w \in F(\bar{p},\bar{x})+\nabla_{x}F(\bar{p},\bar{x})(x-\bar{x})+\partial \tilde{q}(x)\right \}
\end{equation}
around $(0,\bar{x})$. This implication leads immediately to the next statement.
\begin{proposition}\label{Prop.1}
Assume that  $\nabla_{x}F(\bar{p},\bar{x})$ is positive definite. Then $S$ has a single-valued and Lipschitzian localization around $(\bar{p},\bar{x})$.
\end{proposition}

\proof
Observe first that, by \cite[Theorem 13.48]{RW}, $\Sigma$ has the Aubin property around $(0,\bar{x})$ if and only if $(\ref{eq-51})$ has only the trivial solution $u=0$ which, in turn, is ensured by the positive definiteness of $\nabla_{x}F(\bar{p},\bar{x})$. So the assertion follows from  \cite[Theorem 3G.5]{DR} provided the mapping
\[
\Phi(x):= F(\bar{p},\bar{x}) + \nabla_{x}F(\bar{p},\bar{x})(x-\bar{x})+\partial \tilde{q}(x)
\]
is {\em locally monotone} at $(\bar{x},0)$, i.e., for some neighborhood $\mathcal{U}$ of $(\bar{x},-F(\bar{x}))$, one has
\[
\langle x^{\prime}-x, \nabla_{x}F(\bar{p},\bar{x})(x^{\prime}-x) \rangle +
\langle x^{\prime}-x, y^{\prime}-y\rangle \geq 0 ~ \forall ~
(x,y),(x^{\prime}, y^{\prime})\in \gph \partial \tilde{q} \cap \mathcal{U}.
\]
This holds trivially due to the posed assumptions and we are done.
\boxatend
\endproof

In some situations the assumption of positive definiteness of $\nabla_{x}F(\bar{p},\bar{x})$ can be weakened.

\begin{proposition}\label{Prop.2}
Assume that $\tilde{q}$ is convex, piecewise linear-quadratic and the mapping
\begin{equation}\label{eq-53}
\Xi(w):=\left \{ k \in \mathbb{R}^{n} | w \in \nabla_{x}F(\bar{p},\bar{x}) k + \partial \varphi (k)\right \}
\end{equation}
with $\varphi(k):=\frac{1}{2} d^{2}\tilde{q} (\bar{x}| -F(\bar{p},\bar{x}))(k)$ is single-valued on $\mathbb{R}^{s}$. Then $S$ has a single-valued and Lipschitzian localization around $(\bar{p},\bar{x})$.
\end{proposition}

\proof
Again, by virtue of \cite[Theorem 3G.4]{DR} it suffices to
  show  that the single-valuedness of $\Xi$ implies the existence of a single-valued and Lipschitzian localization of $\Sigma$ around $(0,\bar{x})$. Clearly,
\[
\gph \Sigma = \left \{ (w,x) \left | \left [ \begin{split}
& x - \bar{x}\\
& w  - \nabla_{x}F(\bar{p},\bar{x})(x-\bar{x})
\end{split}\right ]
\in \gph \partial \tilde{q} -\left [ \begin{split}
&  \bar{x}\\
&   - F(\bar{p},\bar{x})
\end{split}\right ] \right. \right \}
\]
so that $\Sigma$ is a polyhedral multifunction due to our assumptions imposed on $\tilde{q}$, cf. \cite[Theorem 12.30]{RW}. It follows from \cite{Rob76} (see also \cite[Cor.2.5]{OKZ}) that due to the polyhedrality of $\Sigma$, it suffices to ensure the single-valuedness of $\Sigma(\cdot) \cap \mathcal{V}$ on $\mathcal{U}$, where $\mathcal{U}$ is a convex neighborhood of $0 \in \mathbb{R}^{s}$ and $\mathcal{V}$ is a neighborhood of $\bar{x}$. Let us select these neighborhoods in such a way that
\[
\gph\partial \tilde{q} -
\left [ \begin{split}
& \qquad  \bar{x}\\
&   - F(\bar{p},\bar{x})
\end{split}\right ] = T_{\gph \partial \tilde{q}} (\bar{x}, - F(\bar{p},\bar{x})),
\]
which is possible due to the polyhedrality of $\partial \tilde{q}$. Then one has
\[
\gph \Sigma \cap (\mathcal{U} \times \mathcal{V})=\{(w,\bar{x} +k)\in \mathcal{U} \times \mathcal{V} | w \in \nabla_{x}F(\bar{p},\bar{x})k+D\partial \tilde{q}(\bar{x},- F(\bar{p},\bar{x}))(k)\}.
\]
Under the posed assumptions for any $k \in \mathbb{R}^{n}$
\[
D\partial \tilde{q}(\bar{x},- F(\bar{p},\bar{x}))(k) = \partial \varphi(k),
\]
cf. \cite[Theorem 13.40]{RW}, so that
$\gph \Sigma \cap (\mathcal{U} \times \mathcal{V})=\{(w,k+\bar{x})\in \mathcal{U} \times \mathcal{V} | (w,k)\in \gph \Xi \}$. Since $D \partial \tilde{q}(\bar{x},- F(\bar{p},\bar{x}))(\cdot)$ is positively homogeneous, $\partial \varphi(\cdot)$ is positively homogeneous as well and so the single-valuedness of
$\Sigma(\cdot)\cap \mathcal{V}$ on $\mathcal{U}$ amounts exactly to  the single-valuedness of $\Xi$ on $\mathbb{R}^{s}$.
\boxatend
\endproof
Note that, by virtue of \cite[Theorem 4.4]{GO3}, in the setting of Proposition \ref{Prop.2} one has
\begin{equation}\label{eq-300}
DS(\bar{p},\bar{x})(h)= \{k | 0 \in \nabla_{p}F(\bar{p},\bar{x})h + \nabla_{x}F(\bar{p},\bar{x})k+\partial \varphi (k)\}
\end{equation}
for all $h \in \mathbb{R}^{m}$.

If $\tilde{q} = \delta_{A}$ for a convex polyhedral set $A$, then $\partial \tilde{q}(x)=N_{A}(x)$ and, by virtue of (\ref{eq-102}),               $\varphi(k)=\delta_{\mathcal{K}_{A}(\bar{x},-F(\bar{p},\bar{x}))}(k)$.  It follows that the GE in (\ref{eq-53}) attains the form
\[
w \in \nabla_{x}F(\bar{p},\bar{x})k + N_{\mathcal{K}_{A}(\bar{x},-F(\bar{p},\bar{x}))}(k).
\]
This is in agreement with \cite[Theorem 5.3]{OKZ} and \cite[Theorem 4F.1]{DR}. Let us illustrate the general case of Proposition \ref{Prop.2} via a simple academic\\

\begin{example}
Put $m=2, l=1, n=1$ and consider the GE (\ref{eq-50}), where $F(p,x)=p_{1}+p_{2} x, \tilde{q}(x)=
|x| + \delta_{A}(x)$ with $A=[0,1]$ and the reference pair $(\bar{p},\bar{x})=((-1,1),0)$. Clearly, $\bar{v}:= -F(\bar{p},\bar{x}) = 1$, and we may again employ formula (\ref{eq-102}). One has
 $K(\bar{x},\bar{v})=\mathbb{R}_{+}, f^{\prime \prime} (\bar{x},w)=0$ for any $w\in \mathbb{R}_{+}$ and so we obtain that
\[
\varphi(k)=\frac{1}{2} d^{2} \tilde{q}(\bar{x}|- F(\bar{p},\bar{x}))(k)=\delta_{\mathbb{R}_{+}}(k).
\]
It follows that
\begin{multline}
 \Xi(w)=\{k|w\in k+\partial\delta_{\mathbb{R}_{+}}(k)\}=\{k|w \in k+N_{\mathbb{R}_{+}}(k)\}=  
 \{k \geq 0|k-w\geq 0, \langle k,k-w\rangle =0\}.
\end{multline}

Since $\Xi$ is clearly single-valued on $\mathbb{R}$, we may conclude that the respective $S$ has indeed the single-valued and Lipschitzian localization around $(\bar{p},\bar{x})$.
By virtue of (\ref{eq-300})
\[
DS (\bar{p},\bar{x})(h)=\{k | 0 \in h_{1} + k + N_{\mathbb{R}_{+}}(k)\}.
\]
Both  mappings $S$ and $DS (\bar{p},\bar{x})$ are depicted in Fig.1.
\end{example}

\begin{figure}[h]
\centering
\includegraphics[width=\textwidth]{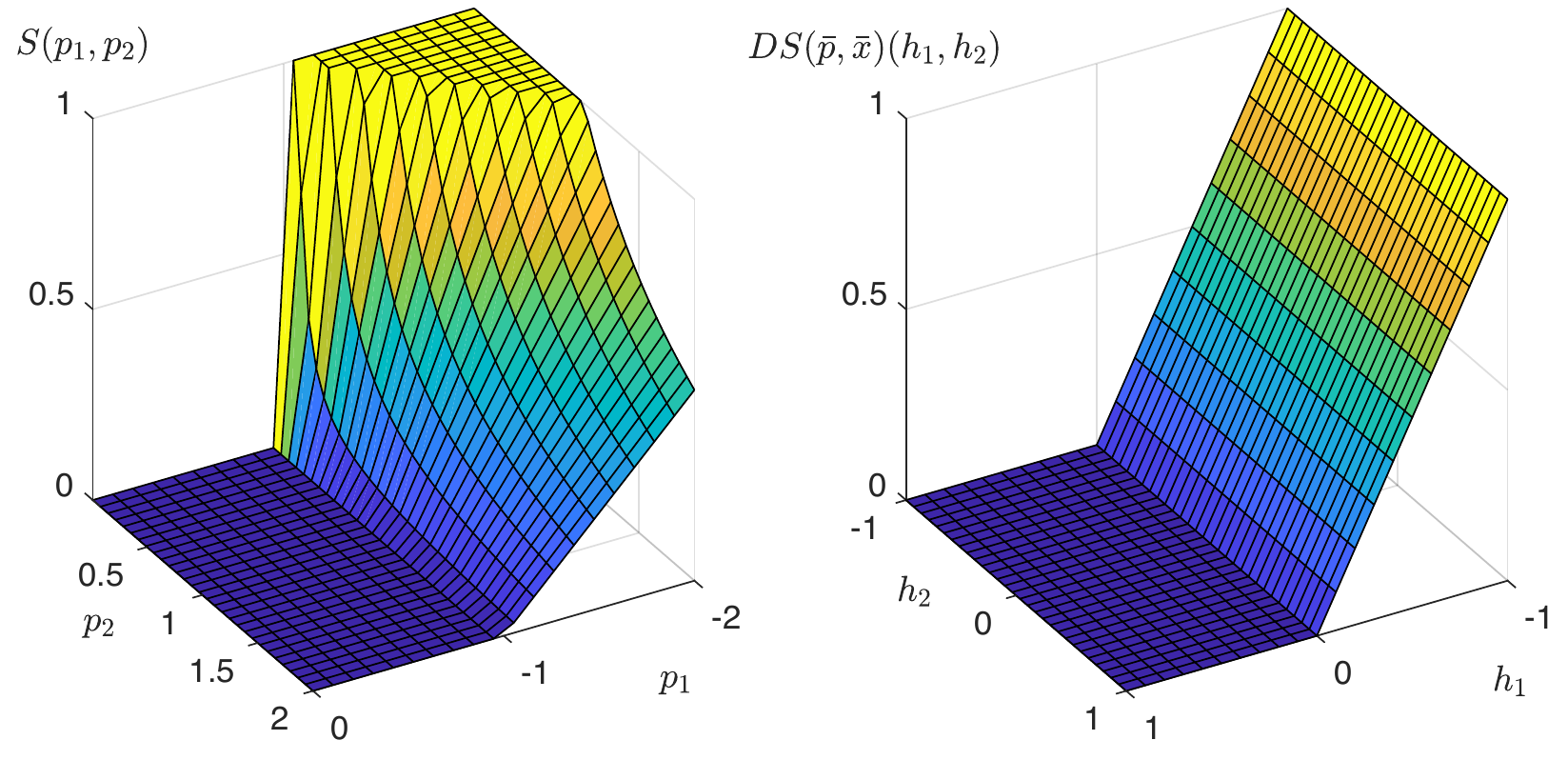}
\caption{The set in equilibria $S(p_1, p_2)$ (left) and the graphical derivative $DS(\bar{p},\bar{x})(h_1,h_2)$  of  Example 1.}
\label{mesh}
\end{figure}

In some cases one can apply the following criterion based on Proposition \ref{Prop.2} and \cite[Theorem 13.9]{RW}.

\begin{proposition}\label{Prop.3}
Assume that $\tilde{q}$ is convex, piecewise linear-quadratic, $ \tilde{q}^{\prime\prime}(\bar{x};\cdot)$ is convex and, with
\[
K:= \{k | q^{\prime}(\bar{x}; k)=\langle -F(\bar{p},\bar{x}), k\rangle \},
\]
one has
\[
\partial \frac{1}{2}d^{2}\tilde{q}(\bar{x}| -F(\bar{p},\bar{x}))(k)= \partial \frac{1}{2} \tilde{q}^{\prime\prime}(\bar{x};k)+N_{K}(k).
\]
Further suppose that the matrix $\nabla_{x}F(\bar{p},\bar{x})$ is copositive with respect to $K$.
Then $S$ has a single-valued and Lipschitzian localization around $(\bar{p},\bar{x})$.

\end{proposition}

\proof
We show that the assumptions of  Proposition \ref{Prop.2} are fulfilled. To this aim consider
the quantity
\[
V:=\langle \nabla_{x} F(\bar{p},\bar{x})(k_{1}-k_{2}), k_{1}-k_{2}\rangle + \langle \xi_{1}-\xi_{2}, k_{1}-k_{2} \rangle + \langle \eta_{1}-\eta_{2}, k_{1}-k_{2}\rangle
\]
with $\xi_{1}\in \partial \frac{1}{2} \tilde{q}^{\prime \prime}(\bar{x}; k_{1}), \xi_{2}\in \partial \frac{1}{2} \tilde{q}^{\prime \prime}(\bar{x}; k_{2}), \eta_{1} \in N_{K}(k_{1}) \mbox{ and } \eta_{2} \in N_{K}(k_{2})$.
It follows that $k_{1}, k_{2} \in K$ so that, by the imposed assumptions, there is a $\sigma > 0$ such that $V  \geq\sigma \| k_{1}-k_{2} \|^{2}$. Hence the operator on the right-hand side of
\[
\xi \in \nabla_{x}F(\bar{p},\bar{x})k + \partial \frac{1}{2} \left[\tilde{q}^{\prime \prime}(\bar{x}; k)\right]+N_{K}(k).
\]
is strongly monotone. This  implies, by virtue of \cite[Proposition 12.54]{RW}, the single-valuedness of the respective mapping $\Xi$ and we are done.  \boxatend
\endproof

 \section{Optimal strategies of producers}
  As stated in the Introduction, our motivation for a study of mapping (\ref{eq-2}) came from an
  attempt to optimize the production strategies of firms with respect to changing external
  parameters like input prices, parameters of inverse demand functions etc. These parameters evolve
  in time and the corresponding adjustments of production strategies have (at least at some
  producers) to take into account the already mentioned costs of change. The appropriate variant of
  the GE in (\ref{eq-2}) (depending on the considered type of market) has thus to be solved at each
  time step with the updated values of the parameters. In this section we will analyze from this
  point of view a standard oligopolistic market  described thoroughly in \cite{MSS} and in \cite[Chapter
  12]{OKZ}. So, in the framework  (\ref{eq-1}) we will assume that $n$ is the number of produced homogeneous commodities,  $p=(p_{1},p_{2})\in \mathbb{R}^{m_{1}} \times
  \mathbb{R}^{m_{2}}, m_{1}+m_{2}=m$ and
  \begin{equation}\label{eq-100}
f_{i} (p,x_{i}, x_{-i})= c_{i}(p_{1},x_{i}) - \langle x_{i}, \pi (p_{2},T)\rangle
  \end{equation}
  with $T=\sum\limits^{l}_{i=1} x_{i}$.
 Functions $c_{i}: \mathbb{R}^{m_{1}} \times \mathbb{R}^{n} \rightarrow
  \mathbb{R}$
represent the {\em production costs} of the $i$th producer and $\pi: \mathbb{R}^{m_{2}} \times \mathbb{R}^{n}
\rightarrow \mathbb{R}$ is the {\em inverse demand function} which assigns each value of the
parameter $p_{2}$ and the overall production vector $T$ the price at which the (price-taking) consumers
are willing to demand. Additionally, we assume that, with some non-negative reals $\beta_{i}$,
  \begin{equation}\label{eq-101}
  q_{i}(x_{i})= \beta_{i}\|x_{i}-a_{i} \|, \quad  i=1,2,\ldots, l,
  \end{equation}
where $||\cdot ||$ stands for an arbitrary norm in $ \mathbb{R}^{n} $.
Sets $A_{i} \subset \mathbb{R}^{n}$ specify the {\em sets of feasible productions} and functions $q_{i}$ represent the
costs of change associated with the change of production from a given vector $a_{i}$ to $x_{i}$.
Thus
$$a_{i}\in A_{i}, \quad i=1,2,\ldots, l,$$ are "previous" productions which have to be changed taking into account the "new" values of parameters $p_{1},p_{2}$.
Clearly, one could definitely work also with more complicated functions $q_i$.
Let us denote the total costs (negative profits) of the single firms by
  \[
  J_{i}(p,x_{i}, x_{-i}):= f_{i}(p,x_{i}, x_{-i})+ q_{i} (x_{i}), \quad i=1,2,\ldots,l.
  \]
In accordance with \cite{MSS} and \cite{OKZ} we will now assume for brevity that $n=1$ (so that $s=l$) and impose
the following assumptions:

\begin{enumerate}
\item [(S1)]
 $\exists$ an open set $\mathcal{B}_{1}\subset \mathbb{R}^{m_{1}}  $ and open sets
 $\mathcal{D}_{i} \supset A_{i}$ such that for for $i=1,2,\ldots, l$
 \begin{itemize}
  \item
  $c_{i}$ are twice continuously differentiable on $\mathcal{B}_{1} \times
  \mathcal{D}_{i}$;
  \item
  $c_{i}(p_{1}, \cdot)$ are convex for all $p_{1} \in \mathcal{B}_{1}$.
 \end{itemize}
  \item [(S2)]
  $\exists$ an open set $\mathcal{B}_{2}\subset \mathbb{R}^{m_{2}}$ such that
  \begin{itemize}
   \item
   $\pi$ is twice continuously differentiable on $\mathcal{B}_{2} \times \inn
   \mathbb{R}_{+}$;
   \item
   $\vartheta \pi (p_{2},\vartheta)$ is a concave function of $\vartheta$ for all $p_{2}\in
   \mathcal{B}_{2}$.
  \end{itemize}
  \item [(S3)]
  Sets $A_{i}\subset \mathbb{R}_{+}$ are closed bounded intervals and at least one of them
  belongs to $\inn \mathbb{R}_{+}$.
\end{enumerate}

Note that thanks to (S3) one has that $T > 0$ for any feasible production profile
$$(x_{1},x_{2},\ldots, x_{l})\in \mathbb{R}^{l} $$ and hence the second term in (\ref{eq-100})
(representing the revenue) is well-defined.

By virtue of \cite[Lemmas 12.1 and 12.2]{OKZ} we conclude that, with $F_{i}$ and $q_{i}$
given by (\ref{eq-100}) and (\ref{eq-101}), respectively, and $\mathcal{B} =  \mathcal{B}_{1}
\times \mathcal{B}_{2}$, the assumptions of Proposition \ref{Prop.1} are fulfilled. This means in
particular that
 for  all vectors $(a_{1},a_{2},\ldots, a_{i})\in A_{1} \times A_{2} \times \ldots \times A_{l}$ the respective mapping $S:(p_{1},p_{2})\mapsto x$  has a single-valued and Lipschitzian localization around any triple $(p_{1},p_{2},x)$, where $(p_{1},p_{2})\in \mathcal{B}_{1}
\times \mathcal{B}_{2}$ and $x \in S(p_{1},p_{2})$. Under the posed assumptions, however, a stronger statement can be established.

\begin{theorem}\label{Thm.2}
Let $a \in A_{1}\times A_{2}\times \ldots \times A_{l}$.
Under the posed assumptions (S1)-(S3) the solution mapping $S$ is single-valued and Lipschitzian over $\mathcal{B}_{1}
\times \mathcal{B}_{2}$.
\end{theorem}

\proof
Given the vectors $a_{i}, i=1,2,\ldots, l$, and the parameters $p_{1}, p_{2}$, the GE in (\ref{eq-2}) attains the form
\begin{equation}
\begin{split}
\label{eq-200}
0 \in
\left[ \begin{array}{c}
\nabla_{x_{1}}c_{1}(p_{1},x_{1})-x_{1}\nabla_{x_{1}}\pi(p_{2},T)-\pi(p_{2},T) \\
\vdots\\
\nabla_{x_{l}}c_{l}(p_{1},x_{l})-x_{l}\nabla_{x_{l}}\pi(p_{2},T)-\pi(p_{2},T)
\end{array}\right]
+  \left[ \begin{array}{c}
\Lambda_{1}(x_{1}-a_{1}) \\
\vdots\\
 \Lambda_{l}(x_{l}-a_{l})
\end{array}\right]
\\ \\
+ ~~ N_{A_{1}}(x_{1})\times \ldots \times N_{A_{l}}(x_{l}),
\end{split}
\end{equation}
where
\[
\Lambda_{i} (x_{i}-a_{i})= \left\{
\begin{split}
\beta_{i} & \mbox{ if }  x_{i} > a_{i}\\
-\beta_{i} & \mbox{ if }  x_{i} < a_{i}\\
[-\beta_{i} & \beta_{i}] \mbox{ otherwise. }
\end{split}
\right.
\]
From \cite[Lemma 12.2]{OKZ} and \cite[Proposition 12.3]{RW} it follows that for any $(p_{1},p_{2}) \in \mathcal{B}_{1}
\times \mathcal{B}_{2}$ the first operator on the right-hand side of (\ref{eq-200}) is strictly monotone in variable $x$. Moreover, the second one, as the subdifferential of a proper convex function is monotone (\cite[Theorem 12.17]{RW}). Their sum is strictly monotone by virtue of \cite[Exercise 12.4(c)]{RW} and so we may recall \cite[Example 12.48]{RW} according to which $S(p_{1}, p_{2})$ can have no more than one element for any  $(p_{1},p_{2}) \in \mathcal{B}_{1}
\times \mathcal{B}_{2}$. This, combined with Theorem \ref{Thm.1} and the Lipschitzian stability of $S$ mentioned above  proves the result.
\boxatend
\endproof

In the next section we will be dealing with the mapping $Z:\mathbb{R}\rightrightarrows     \mathbb{R}^{l-1}$ which, for  given fixed values of $a,p_{1}$ and $p_{2}$, assigns each vector $x_{1} \in A_{1}$ a solution $(x_{2}, \ldots, x_{l})$ of the GE
\begin{equation}\label{eq-400}
\begin{split}
0 \in
\left[ \begin{array}{c}
\nabla_{x_{2}}c_{2}(p_{1},x_{2})-\langle x_{2}\nabla_{x_{2}}\pi(p_{2},T)\rangle -\pi(p_{2},T) \\
\vdots\\
\nabla_{x_{l}}c_{l}(p_{1},x_{l})-\langle x_{l}\nabla_{x_{l}}\pi(p_{2},T)\rangle -\pi(p_{2},T)
\end{array}\right]
+ \left[ \begin{array}{c}
\Lambda_{2}(x_{2}-a_{2}) \\
\vdots\\
 \Lambda_{l}(x_{l}-a_{l})
\end{array}\right]\\
\\
 + ~~  N_{A_{2}}(x_{2})\times \ldots \times N_{A_{l}}(x_{l}).
\end{split}
\end{equation}
Variable $x_{1}$ enters GE (\ref{eq-400}) via $T(=\sum\limits^{l}_{i=1} x_{i})$. Using the same argumentation as in Theorem \ref{Thm.2} we obtain the following result.

\begin{theorem}\label{Thm.3}
Let $a_{i}\in A_{i}$ for $i=1,2,\ldots, l, p_{1}\in \mathcal{B}_{1}$ and $p_{2}\in \mathcal{B}_{2}$. Then, under the assumptions of Theorem \ref{Thm.2}, mapping $Z$ is single-valued and Lipschitzian over $A_{1}$.
\end{theorem}
This statement enables us to consider the situation when the first producer decides to replace the non-cooperative by the Stackelberg strategy, cf. [11, page 220]. In this case, to maximize his profit, he has, for given values of $a, p_{1}$ and $p_{2}$, to solve the MPEC
\begin{equation}\label{eq-401}
\begin{array}{cl}
\mbox{ minimize } & c_{1}(p_{1},x_{1})-\langle x_{1}, \pi(p_{2},T)\rangle + \beta_{1} | x_{1}-a_{1} |\\
x_{1} & \\
\mbox{ subject to } & \\
& x_{-1} = Z(x_{1})\\
& \quad x_{1} \in A_{1}.
\end{array}
\end{equation}

Thanks to Theorem \ref{Thm.3} problem (\ref{eq-401}) can be replaced by the (nonsmooth) minimization problem
 \begin{equation}\label{eq-402}
\begin{array}{ll}
\mbox{ minimize } & \Theta (x_{1})\\
\mbox{ subject to } & \\
& x_{1} \in A_{1}
\end{array}
\end{equation}
in variable $x_{1}$. In \eqref{eq-402}, $\Theta: \mathbb{R}\rightarrow\mathbb{R}$ is the composition defined by
 \begin{equation}\label{eq-403}
\Theta (x_{1}) = c_{l}(p_{1},x_{l})-\langle x_{l}, \pi(p_{2},x_{1}+\mathcal{L}(Z(x_{1})))
\rangle + \beta_{1} | x_{1}-a_{1} |,
\end{equation}
where the mapping $\mathcal{L}:\mathbb{R}^{l-1} \rightarrow\mathbb{R}$ is defined by
\[
\mathcal{L}(x_{2}, x_{3}, \ldots, x_{l}) = \sum\limits^{l}_{i=2} x_{i}.
\]
Problem (\ref{eq-402}) is thus a minimization of a locally Lipschitzian function to which various numerical approaches can be applied.

\section{Numerical experiments}
We consider an example from \cite[Section 12.1]{OKZ}  enhanced by a nonsmooth term reflecting the cost of change.
We have five firms  (i.e., $l=5$) supplying production quantities (productions)
$$ x_1, x_2, \dots, x_5 $$
of one (i.e., $n=1$) homogeneous commodity to a common market with the inverse demand function
$$ \pi(\gamma, T)=5000^{1/ \gamma} T^{-1/ \gamma}, $$
where $\gamma$ is a positive parameter termed {\em demand elasticity}.
In our tests, however, this parameter will be fixed ($\gamma =1 $). The resulting inverse demand function is depicted in Figure \ref{fig_p_and_c} (left).
\begin{figure}
\centering
\includegraphics[width=\textwidth]{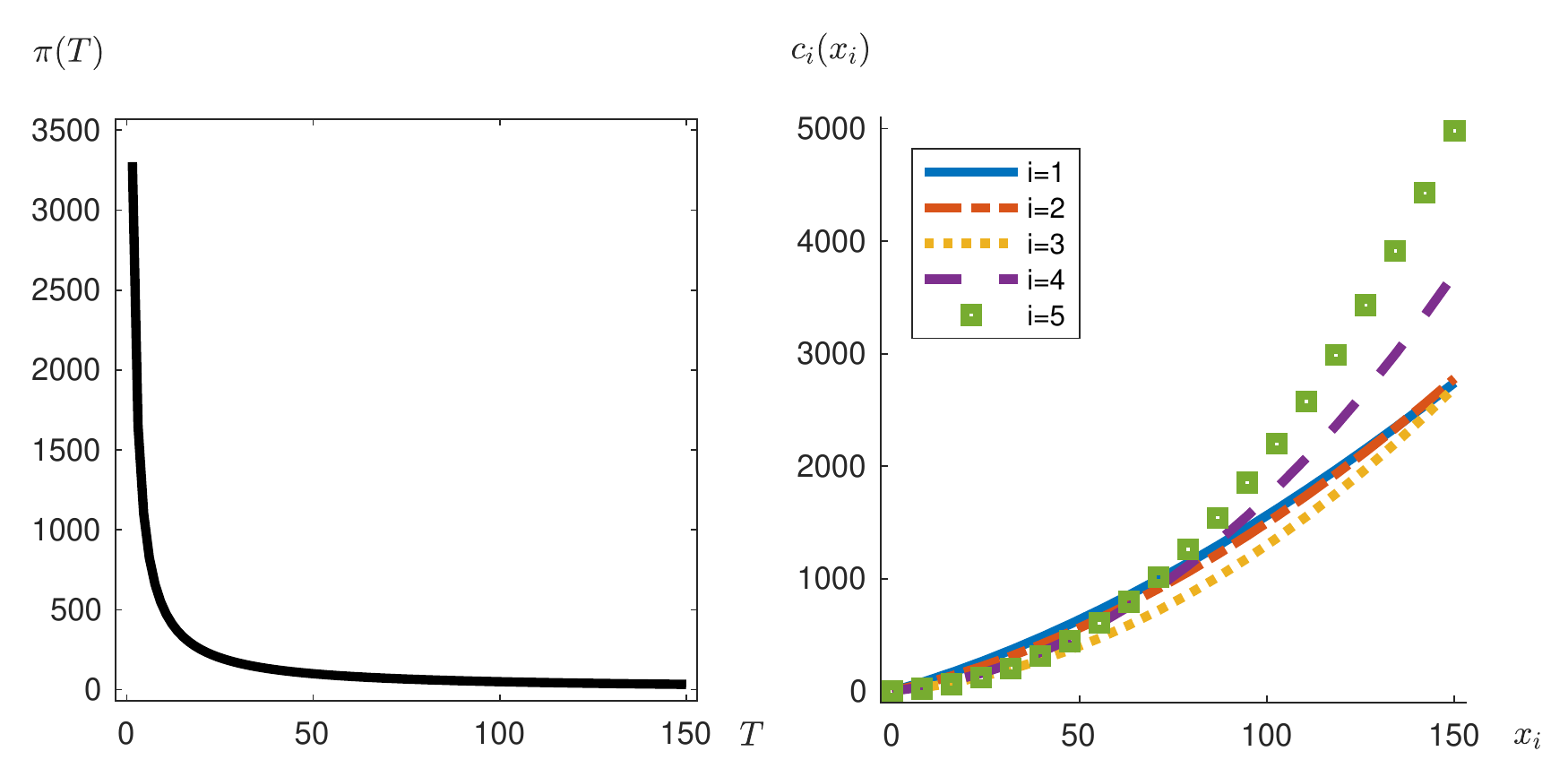}
\caption{The inverse demand function $\pi=\pi(\gamma, T)$  with $\gamma = 1$ (left) and production cost functions $c_i=c_i(x_i), i=1,\dots,5$, for $t=1$ (right).}
\label{fig_p_and_c}
\end{figure}
The production cost functions have the form
$$c_i(b_{i},x_i)=b_i x_i + \frac{\delta_i}{\delta_i+1} {K_i}^{-1/ \delta_i}  x_i^{(1 + \delta_i)/  \delta_i}, $$
where $b_i, \delta_i, K_i, i=1, \dots, 5$, are positive parameters. For brevity we  assume that only the parameters $b_i$, reflecting the impact of the input prices on the production costs, evolve in time, whereas parameters $\delta_i, K_i$ attain the same constant values as in \cite[Table 12.1]{OKZ}.  The cost of change will arise only at the Firms 1, 2 and 3 with different multiplicative constants
$$\beta_1=0.5, \quad \beta_2=1, \quad \beta_3=2.$$
 At the remaining firms any change of production does not incure additional costs ($\beta_4=\beta_5=0$). We will study the behaviour of the market over three time intervals,
 $t \in \{1, 2, 3 \}$
with the initial productions (at $t=0$)
$$a_1=47.81, \quad a_2=51.14, \quad a_3=51.32, \quad a_4=48.55, \quad a_5 = 43.48, $$
corresponding to the standard Cournot-Nash equilibrium with the parameters taken over from \cite{MSS}.

The evolution of parameters $b_i$ is displayed in Table \ref{tab1} and the production cost functions for $t=1$ are depicted in Figure \ref{fig_p_and_c} (right).


\subsection{Cournot-Nash equilibria}

\newcolumntype{b}{X}
\newcolumntype{d}{>{\hsize=1\hsize}X}
\newcolumntype{s}{>{\hsize=.55\hsize}X}

\begin{table}
\normalsize
\begin{tabularx}{\linewidth}{s s d d d s s}
\hline
  & $i$ & 1 & 2 & 3 & 4 & 5 \\
\hline
 t=1 & $b_i$  &  9  & 7  & 3  & 4 & 2 \\
 t=2 & $b_i$ & 10   & 8 & 5  & 4 & 2 \\
 t=3 & $b_i$  & 11  & 9 & 8  & 4 & 2 \\
\hline
\end{tabularx}
\caption{Time dependent input parameters for the production costs.}  \label{tab1}
\normalsize
\begin{tabularx}{\linewidth}{s s d d d s s}
\hline
& $\;i$ & 1 & 2 & 3 & 4 & 5 \\
\hline
t=0 & $\,a_i$ & 47.81& 51.14& 51.32& 48.55& 43.48\\
\hline
t=1 & $\;x_i$ & 49.41& 51.14& 54.24& 48.05& 43.09\\
 & $-J_i$  & 377.23 (-0.80)& 459.95 & 639.95 (-5.83)& 503.44& 507.09\\
\hline
t=2 & $\;x_i$ & 49.41& 51.14& 54.24& 48.05& 43.09\\
 & $-J_i$  & 328.62 & 408.81 & 537.30 & 503.44& 507.09\\
\hline
t=3 & $\;x_i$ & 45.71& 51.14& 51.58& 48.76& 43.64\\
 & $-J_i$  & 286.75 (-1.85)& 379.76 & 386.92 (-5.31)& 527.22& 527.81\\
\hline
\end{tabularx}
\caption{Cournot-Nash equilibria.} \label{tab:Cournot}
\normalsize
\begin{tabularx}{\linewidth}{s s d d d s s}
\hline
& $\;i$ & 1 & 2 & 3 & 4 & 5 \\
\hline
t=0 & $\;a_i$ & 47.81& 51.14& 51.32& 48.55& 43.48\\
\hline
t=1 & $\;x_i$ & 54.95& 51.14& 53.59& 47.52& 42.68\\
 & $-J_i$  & 380.49 (-3.57)& 443.52 & 619.80 (-4.54)& 486.00& 491.88\\
\hline
t=2 & $\;x_i$ & 53.09& 51.14& 53.59& 47.72& 42.84\\
 & $-J_i$  & 329.49 (-0.93)& 398.58 & 523.65 & 492.55& 497.60\\
\hline
t=3 & $\;x_i$ & 53.05& 50.46& 50.77& 48.11& 43.14\\
 & $-J_i$  & 289.65 (-0.02)& 356.57 (-0.68)& 364.33 (-5.64)& 505.29& 508.71\\
\hline
\end{tabularx}
\caption{Stackelberg-Cournot-Nash equilibria. Firm 1 is a leader.} \label{tab:Stackelberg}
\end{table}
For the computation of respective Cournot-Nash equilibria we make use of a ``nonsmooth'' variant of the Gauss-Seidel method described in  \cite[Algorithm 2]{Ka}. Using the notation of Section 4, in the main step S.2 one computes the iteration $x^{k}=(x^{k}_{1}, \ldots, x^{k}_{l})$ via $l$ nonsmooth (but convex) consecutive minimizations
\begin{equation}\label{eq-1000}
x^{k}_{i} \in \argmin\limits_{x_{i}\in A_{i}}f_{i}(p, x^{k}_{1}, \ldots, x^{k}_{i-l}, x_{i},  x^{k-1}_{i+1}, \ldots, x^{k-1}_{l}) 
+q_{i}(x_{i}), i=1,2,\ldots, l,
\end{equation}
where $p$ is the given parameter and $ x^{k-1}=(x^{k-1}_{1}, \ldots,x^{k-1}_{l})$ is the previous iteration.\\
A respective modification of the convergence result \cite[Theorem 6.1]{Ka} takes the following form.
\begin{theorem}\label{Thm.4}
In addition to the assumptions posed in Section 4 suppose that the sequence $x_{k}$
converges for $k \rightarrow \infty$ to some $x^{*}=(x^{*}_{1}, \ldots, x^{*}_{l})\in \mathbb{R}^{l}$. Then $x^{*}$ is a Cournot-Nash equilibrium.
\end{theorem}

\proof
Observe that for $i=1,2,\ldots, l$ the points $x^{k}_{i}$ fulfill the optimality condition
\[
0 \in \nabla_{x_{i}}
f_{i}(p, x^{k}_{1}, \ldots,  x^{k}_{i},   x^{k-1}_{i+l}, \ldots,   x^{k-1}_{l}) + \partial q_{i}(x^{k}_{i})+ N_{A_{i}}(x^{k}_{i}).
\]
Thus, there are elements $\xi^{k}_{i}\in \partial q_{i}(x^{k}_{i})$ and $\eta^{k}_{i}\in N_{A_{i}}(x^{k}_{i})$ such that, by taking subsequences (without relabeling) if necessary, $\xi^{k}_{i} \rightarrow  \xi^{*}_{i},  \eta^{k}_{i} \rightarrow    \eta^{*}_{i}$, satisfying the conditions
\begin{equation}\label{eq-1001}
0 \in \nabla_{x_{i}}
f_{i}(p, x^{*}_{1}, \ldots,  x^{*}_{l} ) + \xi^{*}_{i}+ \eta^{*}_{i}, \xi^{*}_{i} \in \partial q_i(x^{*}_{i}), \eta^{*}_{i} \in N_{A_{i}}(x^{*}_{i}),
\end{equation}
$i=1,2,\ldots, l$. This follows from the outer semicontinuity of the subdifferential mapping and from the uniform boundedness of the subdifferentials $\partial q_{i}$.  Conditions (\ref{eq-1001}) say that $x^{*}$ is a Cournot-Nash equilibrium and we are done.
\boxatend
\endproof

As the stopping rule we employ an approximative version of the optimality conditions (\ref{eq-1001}).

The obtained results are summarized in Table \ref{tab:Cournot} and show  the productions and profits (negative total costs) of all firms at single time instances. In parentheses we display  the costs of change (with negative signs) which decrease the profits of the Firms 1, 2, 3 in case of any change of their production strategy.

Note that firms 1 and 3 increased significantly their productions in time 1 and decreased them in time 3, whereas Firm 2, due to the cost of change, kept its production unchanged during the whole time.

Figure \ref{cournot_optima} shows the total cost functions $J_i$   at time 1. Note that the equilibrium production of Firm 2 lies in a kink point because its cost of change is zero. Expectantly, functions $J_i$ are smooth for $i=4, 5$.\\

\begin{figure}[h]
\centering
\includegraphics[width=0.9\textwidth]{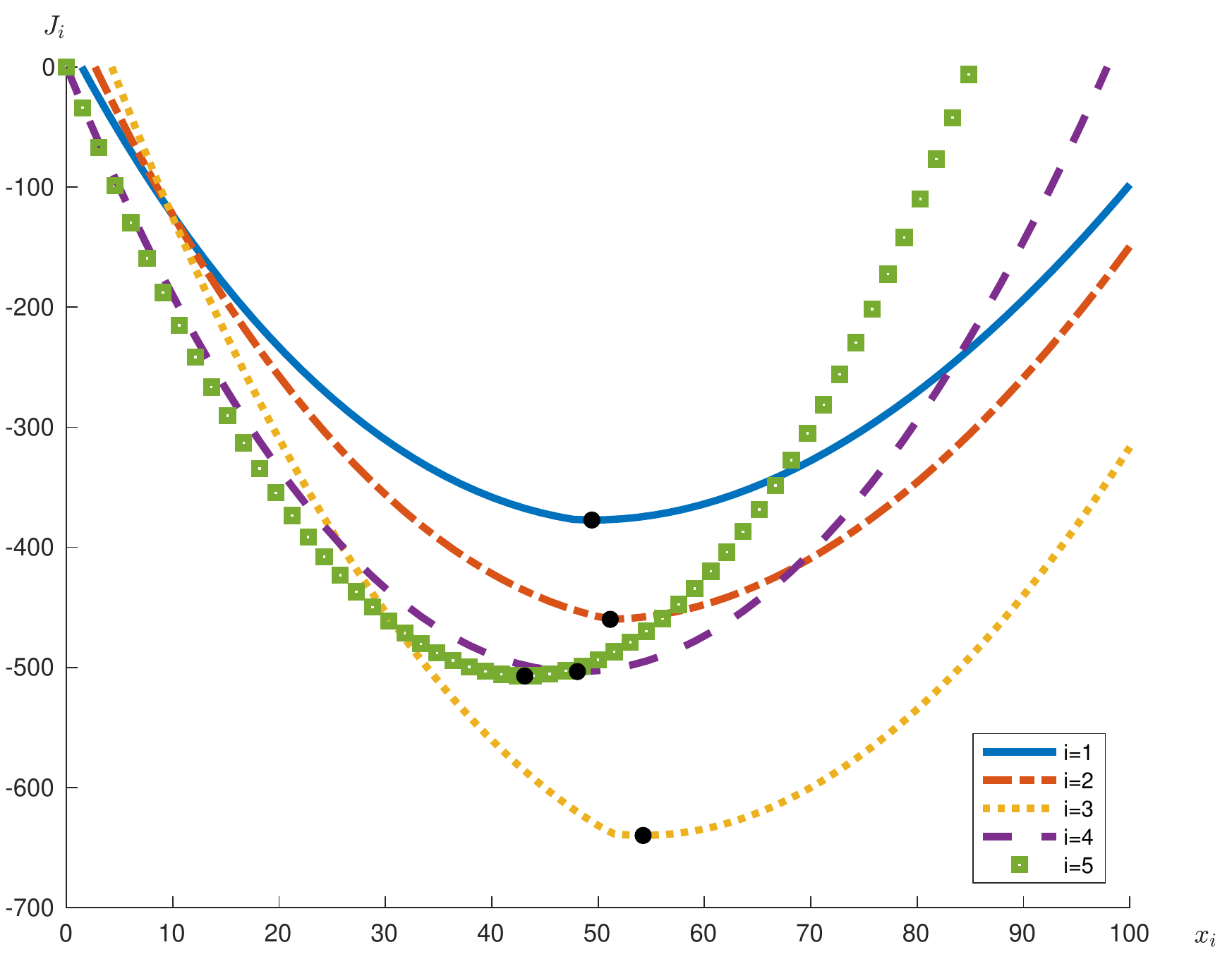}
\caption{The total cost functions and  Cournot-Nash equilibria for t=1. }
\label{cournot_optima}
\end{figure}

\subsection{Stackelberg-Cournot-Nash equilibria}

Next we will consider the same market as in the previous section where, however, the first producer decides now to replace the non-cooperative by the Stackeberg strategy. This leads to the minimization of $\Theta$ over $A_1$,
cf. \cite[page 220]{OKZ}.
 MPEC \eqref{eq-401} is, however, more complicated than the MPECs solved with this approach  in \cite{OKZ} because the presence of costs of change makes the lower-level equilibrium much more difficult and the objective is not continuously differentiable.

On the other hand, since $n=1$, we may use for the minimization of $\Theta$ any suitable routine for nonsmooth constrained univariate minimization and the same routice may be used inside the Gass-Seidel method for the solution of the "lower-level" problem, needed in the computations of the values of $\Theta$.
To this aim we used an  inbuilt Matlab function \verb+fminbnd+ (\cite{Br}).

 \begin{remark}
 Numerical results of both  subsections 5.1, 5.2 were generated by own Matlab code available freely for downloading and testing at:
\begin{center}
{\em \url{https://www.mathworks.com/matlabcentral/fileexchange/72771}} .
\end{center}
The code is flexible and allows for easy modifications to different models.
\end{remark}

The obtained results are summarized in Table 3. They are quite different from their counterparts in Table 2 and show that, switching to the Stackelberg strategy, Firm 1 substantially improves its profit. In contrast to the noncooperative strategy, it has now to change its production at each time step and also Firm 2, who preserved in the Cournot-Nash case the same production over the whole time, is now forced to change it at $t=2$. Of course, our data are purely academic and can hardly be used for some economic interpretations. On the other hand, the results are sound and show the potential of the suggested techniques to be applied also in some more realistic situations.\\


\section*{Conclusion}

\noindent In the first half of the paper we have studied a parametrized variational inequality of the second kind. In this form,  one can write down, for example, a condition which characterizes solutions of some parameter-dependent Nash equilibrium problems. By using standard tools of  variational analysis sufficient conditions have been derived ensuring the existence of a single-valued and Lipschitzian localization of the respective solution mapping. Apart from post-optimal analysis, the obtained results can be used in computation of respective equilibria for given values of the parameter via continuation (\cite{AG}) or Newton-type methods (\cite{GO8}).

The second part of the paper has been inspired, on one hand, by the successful  theory of rate-independence processes (\cite{MR},\cite{FKV}) and, on the other hand, by the important economic paper \cite{Fl}. It turns out that in some market models the cost of  change of the production strategy can be viewed as the economic counterpart of the dissipation energy, arising in rate independent dissipative models of nonlinear mechanics of solids. Cost of change (dissipation energy) occurs further, e.g., in modeling the behavior of some national banks who try to regulate the inflation rate, among other instruments, via buying or selling suitable amounts of the domestic currency on international financial markets \cite{R}.

The considered market model,  obtained by augmenting the cost of change to a standard model from \cite{MSS} possesses very strong stability properties and is amenable to various numerical approaches. Both in the case of the Cournot-Nash equilibrium as well as in the two-level (Stackelberg) case we have used rather simple techniques based on the "nonsmooth" line search from \cite{Br}. In case of multiple commodities and/or more complicated costs of change one could employ more sophisticated approaches based either on the second-order subdifferentials or on the second subderivatives discussed in Section 3. This will be the subject of our future work.

\end{document}